\documentclass[11pt]{amsart}
\usepackage{amssymb}
\usepackage{amsmath}
\usepackage{amsfonts}
\newtheorem{theorem}{Theorem}
\theoremstyle{plain}
 \newcommand{\fin}{\hfill  $\square$}

\newtheorem {lemma}[theorem]{Lemma}
\newtheorem {prop}[theorem]{Proposition}

\newtheorem {remarks}[theorem]{Remarks}

\newtheorem {cor}[theorem]{Corollary}
\begin{document}

\title[ ]{Complex moment problem and recursive relations}

\author{ K. Idrissi }
\address{  K. Idrissi, Mohamed V University. Rabat Morocco }
\email{kaissar.idrissi@gmail.com }

\author{E. H. Zerouali}
\address{E. H. Zerouali, Mohamed V University. Rabat Morocco}
\email{zerouali@fsr.ac.ma}
\date{}

\begin{abstract}
 We introduce a new strategy in  solving the truncated complex moment problem. To this aim we investigate   recursive doubly indexed  sequences and their characteristic polynomials. A characterization of recursive doubly indexed   \emph{moment} sequences is provided. As a   simple application, we obtain  a computable solution to the complex moment problem for cubic harmonic characteristic polynomials of the form $z^3+az+b\overline{z}$, where $a$ and $b$ are arbitrary real numbers. We also recapture a recent result due to Curto-Yoo given for cubic column relations in $M(3)$ of the form $Z^3=itZ+u\overline{Z}$ with  $t,u$  real numbers satisfying some suitable inequalities.  Furthermore, we solve the truncated complex moment problem with column dependence relations of the form $Z^{k+1}= \sum\limits_{0\leq n+ m \leq k} a_{nm} \overline{Z}^n Z^m$ ($a_{nm} \in \mathbb{C}$).
  \end{abstract}

\maketitle

\section{ Introduction}
Let   $\gamma = \{\gamma_{ij}\}_{i,j\geq 0}$  be a doubly  indexed  sequence of complex numbers such that  $\overline{\gamma}_{ij}=\gamma_{ji} $ and $\gamma_{00}>0$. The truncated complex moment problem (TCMP for short)  associated with  $\{\gamma_{ij}\}_{0\leq i, j\leq r}$  entails finding a positive Borel measure $\mu$ supported in the complex plane $\mathbb{C}$ such that
 \begin{equation}\label{i.0}
\gamma_{ij}=\int \overline{z}^iz^jd\mu \phantom{du texte} (0\leq i, j\leq r);
\end{equation}
A sequence  $\{\gamma_{ij}\}_{0\leq i, j\leq r}$   satisfying (\ref{i.0})   will be  called a truncated moment sequence and $\mu$ is said to be  a representing measure for $\{\gamma_{ij}\}_{0\leq i, j\leq r}$.
The full complex moment problem (FCMP) prescribes moments of all orders.
 More precisely, an infinite doubly sequence   $\{\gamma_{ij}\}_{i,j\geq 0}$, with $i, j \in \mathbb{Z}_+$, is a moment sequence provided that there exists a  Borel measure $\mu$ supported in the complex plane $\mathbb{C}$ such that,
 \begin{equation}\label{i.f}
\gamma_{ij}=\int \overline{z}^iz^jd\mu  \text{ for all }i,j\in\mathbb{Z}_+.
\end{equation}
  In \cite{1,2,3} Curto-Fialkow introduced an approach to study the existence and uniqueness of solutions of the TCMP, $\gamma^{(2n)} :=\{\gamma_{ij}\}_{0\leq i+j\leq 2n}$ with $\overline{\gamma}_{ij}=\gamma_{ji}$ and $\gamma_{00}>0$,   based on positivity and extensions  of the moment matrix $M(n)\equiv M(n)(\gamma^{(2n)})$,  built as follows.

\begin{equation*}
M(n)=\left( \begin{matrix}
 M[0,0]&M[0,1]&\ldots&M[0,n]\\
 M[1,0]&M[1,1]&\ldots&M[1,n]\\
 \vdots&\vdots&\ddots&\vdots\\
 M[n,0]&M[n,1]&\ldots&M[n,n]
 \end{matrix} \right),
\end{equation*}
where
\begin{equation*}
M[i,j]=\left( \begin{matrix}
 \gamma_{i,j}&\gamma_{i+1,j-1}&\ldots&\gamma_{i+j,0}\\
 \gamma_{i-1,j+1}&\gamma_{i,j}&\ldots&\gamma_{i+j-1,1}\\
 \vdots&\vdots&\ddots&\vdots\\
 \gamma_{0,i+j}&\gamma_{1,i+j-1}&\ldots&\gamma_{j,i}
 \end{matrix} \right).
\end{equation*}
The matrix $M(n)$ detects the positivity of the Riesz functional
$$
\Lambda_{\gamma^{(2n)}}:p(Z,\overline{Z}) \equiv \sum_{0\leq i+j\leq 2n}a_{ij}\overline{Z}^iZ^j \longrightarrow \sum_{0\leq i+j\leq 2n} a_{ij}\gamma_{ij}
$$
 on the cone generated by the collection $\{p\overline{p}: p\in\mathbb{C}_n[Z,\overline{Z}]\}$, where $\mathbb{C}_n[Z,\overline{Z}]$ is the vector space of polynomials in two variables with complex coefficients and total  degree less than or equal to $n$. In the sequel, we will write$,  d_z(P), d_{\bar z}(P)$ and $d_P \equiv \deg P$ for the degree in $z$, the degree in $\bar z$ and the total degree of $P$ respectively.
  Considering the lexicographic order
$$1, Z, \overline{Z}, Z^2, Z\overline{Z}, \overline{Z}^2, \dots, Z^n ,Z^{n-1}\overline{Z}, \ldots ,Z\overline{Z}^{n-1}, \overline{Z}^n$$
to denote rows and columns of the moment matrix $M(n)$. It is immediate  that the rows $\overline{Z}^kZ^l$, columns $\overline{Z}^iZ^j$ entry of the matrix $M(n)$ is equal to $\Lambda_{\gamma^{(2n)}}(\overline{z}^{i+l}z^{j+k}) =\gamma_{i+l,j+k}$. For reason of simplicity, we identify a polynomial $p(z; \overline{z})\equiv \sum a_{ij} \overline{z}^i z^j$ with its coefficient vector $p = (a_{ij})$ with respect to the basis of monomials of $\mathbb{C}_n[z; \overline{z}]$ in degree-lexicographic order. Clearly  $M(n)$ acts on these coefficient vectors as follows: $$<M(n) p, q> = \Lambda_{\gamma^{(2n)}}(p\overline{q}).$$

A result of Curto-Fialkow \cite{1}  states that $\gamma^{(2n)}=\{\gamma_{ij}\}_{0\leq i+j\leq 2n}$ has a representing measure if and only if $M(n)$ admits a positive finite rank moment matrix extension $M(\gamma) \equiv M(\infty)$. In general the existence of such extension is difficult to determine, but in
various special cases concrete results are known. In particular, a complete solution to TCMP based on moment matrix extensions is known for $n \leq 2$ and for $M(n)$ whenever the submatrix $M(2)$ is singular (cf. \cite{cy1}, \cite{cr1}, \cite{5}). Therefore, the first open case of TCMP concerns $n = 3$ with $M(2)$ positive definite. The Curto-Yoo paper \cite{4} concerns part of this problem. In general, if $\gamma^{(2n)}$ has a representing measure, then $M(n)$ is positive and $\text{rank } M(n) \leq \text{card } \nu_{\gamma^{(2n)}}$, where $\nu_{\gamma^{(2n)}}$, the variety associated with $\gamma^{(2n)}$  as the intersection of the zero-sets of the polynomials $p(z;\overline{z})$ such that $M(n) p= 0$. In Curto-Fialkow-Moller \cite{cfm}, TCMP is solved for the "extremal" case when $\text{rank }M(n) = \text{card }\nu_{\gamma^{(2n)}}$. In this case, $\gamma^{(2n)}$ has a representing measure if and only if $M(n)\geq 0$ and $\gamma^{(2n)}$ is consistent. That is, if $p \in \mathbb{C}_{2n}[z; \overline{z}]$ and $p_{\mid \nu_{\gamma^{(2n)}}} \equiv 0$, then $\Lambda_{\gamma^{(2n)}}(p) = 0$. The proof of this fact  does not require results on moment matrix extensions; it is based  on elementary  tools  from  vector space duality, convex analysis, and interpolation by polynomials.

 As noted above, the simplest unsolved case of TCMP concerns $M(3)(\gamma^{(6)})$. Within this problem, \cite{4} identifies a subproblem which is "extremal", with $\text{rank } M(3) = \text{card }\nu_{\gamma^{(6)}} = 7$. Indeed, \cite{4} focuses on the case when $M(3)\geq 0$ and $M(2)>0$, and there is a column dependence relation $p(Z;\overline{Z}) = 0$ (and automatically, $\overline{p}(Z;\overline{Z}) = 0$), where $p$ is a harmonic polynomial of the form $p(z,\overline{z})= z^3-itz-u\overline{z}$, with real parameters $t$ and $u$. The $M(3)$ problem remains open and is only partially solved. Lemma 2.3 in \cite{4} states that if $0 < \mid u\mid < t < 2\mid u \mid$, then $p$ has exactly 7 zeros in the complex plane. This later can be disproved by  the next two examples.\\

 \noindent  {\bf Example 1.} The equality $z^3=2iz-\frac{5}{4}\overline{z}$ admits only 3 zeros, 0 and $\pm\frac{\sqrt{13}}{2}e^{i\frac{\pi}{4}}$, although $t=2$ and $u=-\frac{5}{4}$ verify the condition $0<|u|<t<2|u|$.\\

 \noindent  {\bf Example 2.} The equality $z^3=-2iz+\frac{5}{4}\overline{z}$ admits 7 zeros, 0 ; $\pm\frac{\sqrt{3}}{2}e^{-i\frac{\pi}{4}}$ and $(\pm\frac{3\sqrt{2}}{4}\pm i\frac{\sqrt{2}}{4})e^{-i\frac{\pi}{4}}$, although $t=-2$ and $u=\frac{5}{4}$ does not verify the condition $0<|u|<t<2|u|$.\\

 This "inattention" led to the next  incorrect version of the main theorem in \cite{4},

 \begin{theorem}\label{c.m.t}\cite[Theorem 1.5]{4}

Let $M(3)\geq0$, with $M(2)> 0$ and $z^3-itz-u\overline{z}=0$. For $u,t\in\mathbb{R}$, assume that $0<|u|<t<2|u|$. The following statements are equivalent.
\begin{itemize}
  \item [$i)$] There exists a representing measure for $M(3)$.

  \item [$ii)$]
    \begin{equation*}\begin{cases}
\Lambda_{\gamma^{(6)}}(q_{LC})&=0,\\
\Lambda_{\gamma^{(6)}}(zq_{LC})&=0.
    \end{cases}\end{equation*}

  \item [$iii)$]
    \begin{equation*}\begin{cases}
\mathfrak{Re}\gamma_{12}-\mathfrak{Im}\gamma_{12}&=u(\mathfrak{Re}\gamma_{01}-\mathfrak{Im}\gamma_{01}),\\
\gamma_{22}&=(t+u)\gamma_{11}-2u\mathfrak{Im}\gamma_{02}.
    \end{cases}\end{equation*}

  \item [$iv)$] $q_{LC}:=i(z-i\overline{z})(z\overline{z}-u)=0$,
\end{itemize}
\end{theorem}

 Example 1 shows that for $u < 0$ the last theorem is not valid, because $6 \leq \text{rang }M(3) $ and $\text{card } \nu_{\gamma^{(6)}} \leq 3$ ( recall that, if $\gamma^{(6)}$ is a moment sequence then $\text{rang }M(3) \leq \text{card } \nu_{\gamma^{(6)}}$ ). Theorem  \ref{c.m.t} implies  Corollary 4.3  in the same paper, which has inherited the same mistakes.

 From the previous discussion, it appears natural     to give a new solution to  the TCMP for cubic column relations in $M(3)$ of the forms $Z^3 = itZ + u\overline{Z}$ and $Z^3 = a Z + b \overline{Z}$, where $a,b,t$ and $u$ are real numbers. This is the main goal of  Section 5,  using our methodology on recursive sequences.\\

  We notice below  that the recursiveness in the truncated moment problem is a natural concept and is totally  inherent.   It is obvious  that the truncated moment problem is equivalent to the  recursive full moment problem. Indeed, given a doubly indexed truncated moment sequence $\omega \equiv \{ \gamma_{ij}\}_{ 0\leq i, j\leq n}$. A result of C. Bayer and J. Teichmann \cite{bt} states  that if a finite double sequence of complex numbers has a representing measure, then it has a finitely atomic representing measure. It follows that $\omega$ admits a finite support representing measure $\mu$, suppose that $\text{supp}(\mu) \subset \nu= \{\lambda_1,\lambda_2, \ldots,\lambda_r\} \subset \mathbb{C}$. Let $p\equiv z^{s+1} - \sum\limits_{i+j=0}^{s} a_{ij}\overline{z}^i z^j$ be a polynomial vanishing on $\nu$, since $\int \overline{z}^n z^m p(z, \overline{z}) d\mu =0$ for all $n, m \in \mathbb{Z}_+$, then $\mu$ is a representing measure for the recursively generated  sequence $\gamma\equiv\{\gamma_{ij}\}_{i,j\geq0}$ ($\omega \subset \gamma$) defined by the next relations,

 \begin{itemize}
\item [$i)$] For all $i,j\geq0$,
 \begin{equation}\label{a.1}
 \gamma_{ji}=\overline{\gamma}_{ij}.
 \end{equation}
\item [$ii)$]For all $n,m\geq0$,
 \begin{equation}\label{a.2}
 \gamma_{n, m+s+1}= \sum\limits_{i+j=0}^{s} a_{ij}\gamma_{n+i,m+j}.
 \end{equation}
\end{itemize}
 We shall refer to the polynomial $p(z,\overline{z})= z^{s+1} - \sum\limits_{i+j=0}^{s} a_{ij}\overline{z}^i z^j$ as a  characteristic polynomial associated with $\gamma$. The sequence $\gamma$ can be associated with several characteristic polynomials, see Section 2.

 Thus   every truncated complex moment sequence is a subsequence of a recursively moment sequence.  We deduce that  solving the TCMP is actually equivalent to solve the recursive  full  moment problem. The main goal in this paper is to investigate use the recursive double sequences (verifying \eqref{a.1} and \eqref{a.2}) to get an approach, based on the localization of zeros of the characteristic polynomials, for solving the TCMP.\\

 This paper is organized as follows. We define in Section 2  recursive doubly indexed  sequences. We show that for such sequences the TCMP and FCMP are equivalent. We devote Section 3 to the study of  RDIS associated with analytic characteristic polynomial and we show that the family of analytic characteristic polynomials is a principal ideal of $\mathbb{C}[X]$. Section 4 is devoted to give a characterization of recursive   doubly indexed  moment sequences. In the last section, we apply  our results to give an explicit  solution for the TCMP with cubic column relations in $M(3)$ of the form $Z^3=aZ+b\overline{Z}$ ($a,b\in \mathbb{R}$) and also to regain the correct form of Theorem \ref{c.m.t}. In the last section, we involve the RDIS to give a necessary and sufficient condition for the existence of a solution to the TCMP with column dependence relations of the form $Z^{k+1}= \sum\limits_{0\leq n+ m \leq k} a_{nm} \overline{Z}^n Z^m$ ($a_{nm} \in \mathbb{C}$).

\section{ Recursive Double indexed  sequences}

 Let $\{a_{lk}\}_{0\leq l,k\leq r}$ be some fixed complex numbers and let $\gamma\equiv\{\gamma_{ij}\}_{i,j\geq0}$ be a doubly indexed
 sequence, with $\gamma_{00}>0$, defined by the following relations:

\begin{itemize}
\item [$i)$] For all $i,j\geq0$,
 \begin{equation}\label{1.1}
 \gamma_{ji}=\overline{\gamma}_{ij}.
 \end{equation}
\item [$ii)$] For all $i$ and $n$ with $0\leq i$ and $r\leq n$,
 \begin{equation}\label{1.2}
 \gamma_{i,n+1}=\sum\limits_{0\leq l+k\leq r} a_{lk} \gamma_{l+i,n+k-r}.
 \end{equation}
\end{itemize}
 Where $\omega \equiv \{\gamma_{ij}\}_{0\leq i\leq j\leq r}$ are given initial conditions.

 In the sequel we shall refer to  such  sequence as Recursive Double Indexed  Sequence, RDIS for short. The polynomial
 $P(z,\overline{z})= z^{r+1}-\sum\limits_{0\leq l+k\leq r} a_{lk}\overline{z}^lz^{k}$ is called a characteristic polynomial associated with $\gamma$,
  given by \eqref{1.1} and \eqref{1.2}. This last polynomial has a finite roots set. More precisely $P$ has at most $(r+1)^2$ roots
  ( see Proposition 4.4 in \cite{9}).

 A  RDIS can be defined in various ways using different characteristic polynomials as is shown in the following example.
 Let $\gamma=\{\gamma_{ij}\}_{i,j\geq0}$ with $\gamma_{ij}= \frac{(-1)^{i+j}}{2}+ \frac{1}{2} \mathfrak{Re}((1-2i)^i(1+2i)^j)$.
 Then $\gamma$ may be defined by the following recursive relations,

 \begin{itemize}
   \item $\gamma_{n+2,m}= -2\gamma_{n,m+1}- \gamma_{n,m}$, with $\gamma_{n,m}=\overline{\gamma}_{m,n}$,
    for $\gamma_{00}=1, \gamma_{01}= \gamma_{10}= 0, \gamma_{11}= 3$.
   \item $\gamma_{n+3,m}= \gamma_{n+2,m}-3\gamma_{n+1,m}-5\gamma_{n,m}$, with $\gamma_{n,m}=\overline{\gamma}_{m,n}$,
    for $\gamma_{00}=1, \gamma_{01}= \gamma_{10}= 0, \gamma_{11}= 3, \gamma_{02}= \gamma_{20}= -1, \gamma_{21}= \gamma_{12}= 2, \gamma_{22}= 13.$
 \end{itemize}

 Therefore, $P_1(z,\overline{z})= z^2+2\overline{z}+1 $ and $P_2(z,\overline{z})= Q(z)= z^3-z^2+3z+5$ are two characteristic polynomials associated with $\gamma$.

 Let $\mathcal{P}_\gamma$ be the set of characteristic polynomials associated with the sequence $\gamma\equiv\{\gamma_{ij}\}_{i,j\geq0}$.

\begin{remarks}\label{rmk}
\begin{itemize}
\item [$i)$] The subset $\mathcal{P}_\gamma$ is an ideal of $\mathbb{C}[z,\overline{z}]$.

\item [$ii)$]  The characteristic polynomial $P$, together with the initial conditions and the relations \eqref{1.1} and \eqref{1.2},
   are  said to define the sequence $\gamma$.

\item [$iii)$]  We notice that because of condition \eqref{1.1} ,  Equation  \eqref{1.2} is equivalent to : For all $n$ and $j$
    with $0\leq j$ and $r\leq n$,
 \begin{equation}\label{1.3}
 \gamma_{n+1,j}=\sum\limits_{0\leq l+k\leq r} \overline{a}_{lk} \gamma_{n-r+k,l+j}.
 \end{equation}
 The polynomial $Q(z,\overline{z})=\overline{z}^{r+1}-\sum\limits_{0\leq l+k\leq r} \overline{a}_{lk}\overline{z}^{k}z^l $ is a characteristic polynomial
 associated with $\{\gamma_{i,j}\}_{i,j\geq0}$ given by \eqref{1.3}, where $\{\gamma_{ij}\}_{0\leq j\leq i\leq r}$ are given initial conditions.\\
\end{itemize}
\end{remarks}

 The following result is an immediate consequence of \eqref{1.2}.
 \begin{lemma}\label{l.0}
 Let $\gamma\equiv\{\gamma_{ij}\}_{i,j\geq0}$ be a doubly indexed sequence and let $p(z,\overline{z}) \in \mathbb{C}[z, \overline{z}]$. Then $p(z,\overline{z})$ is a characteristic of $\gamma$ if and only if $M(\gamma) p =0$.
 \end{lemma}

 We use a structural property of moment matrices to get the following interesting result.

\begin{lemma} Under the notations above, for every  $f, g, h \in \mathbb{R}[x_1, \ldots, x_d]$, we have
  \begin{equation}\label{p3.1}
  f^T M(\gamma) (gh) = (fg)^T M(\gamma) h.
  \end{equation}
\end{lemma}

  {\it Proof.} Let  $f, g, h \in \mathbb{R}[x_1, \ldots, x_d]$ be polynomials. We write $f =\sum\limits_{\mathbf{i}} f_{\mathbf{i}} \mathbf{x}^{\mathbf{i}}$, $g =\sum\limits_{\mathbf{j}} g_{\mathbf{j}} \mathbf{x}^{\mathbf{j}}$ and
  $h =\sum\limits_{\mathbf{k}} h_{\mathbf{k}} \mathbf{x}^{\mathbf{k}}$. As the entry of the moment matrix corresponding to the column $\mathbf{x}^{\mathbf{i}}$ and the line $\mathbf{x}^{\mathbf{j}}$ is $\gamma_{\mathbf{i}+\mathbf{j}}$, we obtain
  $$\begin{array}{lll}
f^T M(\gamma) (gh)
&  = & (\sum\limits_{\mathbf{i}} f_{\mathbf{i}} \mathbf{x}^{\mathbf{i}})^T M(\gamma) (\sum\limits_{\mathbf{j}, \mathbf{k}} g_{\mathbf{j}} h_{\mathbf{k}}
\mathbf{x}^{\mathbf{j}+\mathbf{k}})\\
  &=& \sum\limits_{\mathbf{i}, \mathbf{j}, \mathbf{k}} f_{\mathbf{i}} g_{\mathbf{j}} h_{\mathbf{k}} \gamma_{ \mathbf{i}+\mathbf{j}+\mathbf{k} }\\
  &=& (fg)^T M(\gamma) h. \end{array}$$
It follows that
\begin{prop}\label{l.4}
  Let $\gamma\equiv\{\gamma_{ij}\}_{i,j\geq0}$ be a doubly indexed sequence and let $M(\gamma)$ such that its associated moment matrix $M(\gamma)$ is semidefinite positive. Then, for any polynomial $p \in \mathbb{C}[z, \overline{z}]$ and any integer $n\geq 1$,
  \begin{equation}\label{p3.2}
  M(\gamma) p^n =0  \Longrightarrow   M(\gamma) p =0.
  \end{equation}
  \end{prop}

 {\it Proof.}    If $M(\gamma) p^2 =0$, then $0 = M(\gamma) p^2 = 1^T M(\gamma) p^2 = p^T M(\gamma) p$, from \eqref{p3.1}; since $M(\gamma) \geq 0$,  we obtain $M(\gamma) p = 0$ and hence \eqref{p3.2} holds for $n=2$. By induction, \eqref{p3.2} remains valid for any power of $2$. Now , if $M(\gamma) p^n =0$ we choose $r$ in such a way that $r+k$ is a power of $2$ to ensure that
  $$M(\gamma) p^{n+r}=( p^r)^\perp M(\gamma) p^n =0.$$
Which gives $M(\gamma) p =0$.\\

  In the next, we involve the celebrated Hilbert's Nullstellensatz to obtain a very useful result. Let $I_p \equiv (p)$ be the ideal of $\mathbb{C}[z,\overline{z}]$ generated by $p$. The set $V(I_p):=\{z\in \mathbb{C} \mid f(z)=0 \text{ for every } z\in I_p\}$ is the (complex) variety associated with  $I_p$.   $I(V(I_p)):=\{f(z) \in \mathbb{C}[z,\overline{z}] \mid f(x)=0 \text{ for every } z\in V(I_p) \}$ and
 $\sqrt{I_p} := \{ f \in \mathbb{C}[z,\overline{z}]  \mid f^k \in I_p \text{ for some integer} k \geq 1 \}$, are again ideals in $\mathbb{C}[z,\overline{z}]$,  that  obviously contain the ideal $I_p$. The Hilbert's Nullstellensatz states that $$I(V(I_p)) = \sqrt{I_p}.$$

 Now let us consider a polynomial
 $q \in \mathbb{C}[z,\overline{z}]$ satisfying that $Z(p):=\{z\in {\mathbb C}~~\text{such that}~~p(z,\overline{z})=0\} \subseteq Z(q)$. We have  $q \in I(V(I_p))$ hence $q \in \sqrt{I_p}$, that is, there exists some integer
 $k \geq 1$ such that $q^k \in I_p$. Thus, by Remark \ref{rmk}-$i)$ and Lemma \ref{l.0}, $M(\gamma) q^k =0$. This implies, by Proposition \ref{l.4}, that $M(\gamma)q =0$ and therefore, from Lemma \ref{l.0}, $q$ is a characteristic polynomial of $\gamma$.

 \begin{prop}\label{p.11}
 Let $\gamma\equiv\{\gamma_{ij}\}_{i,j\geq0}$ be a RDIS, associated with the characteristic polynomial $p(z, \overline{z})$. Then every polynomial
 vanishing at all points of $Z(p):=\{z\in {\mathbb C}~~\text{such that}~~p(z,\overline{z})=0\}$ is a characteristic polynomial of $\gamma$.
 \end{prop}

 In the following proposition, we show that for a given RDIS the TCMP and the FCMP are equivalent.

\begin{prop}\label{p.1}
 Let $\gamma\equiv\{\gamma_{ij}\}_{i,j\geq0}$ be a RDIS whose initial conditions and characteristic polynomial are $\{\gamma_{ij}\}_{0\leq i\leq j\leq r}$
 and $P(z,\overline{z})=z^{r+1}-\sum\limits_{0\leq l+k\leq r} a_{lk}\overline{z}^lz^{k}$, respectively. The following are equivalent.
 \begin{itemize}
  \item [$i)$]There exists a positive Borel measure $\mu$, solution of the FCMP for $\gamma\equiv\{\gamma_{ij}\}_{i,j\geq0}$.
  \item [$ii)$] There exists a positive Borel measure $\mu$, solution of the TCMP for $\omega\equiv\{\gamma_{ij}\}_{0\leq i\leq j\leq r}$ with
  $$\text{supp}(\mu)\subset Z(P):=\{z\in {\mathbb C}~~\text{such that}~~P(z,\overline{z})=0\}.$$
 \end{itemize}
\end{prop}

{\it Proof.}
 \begin{itemize}
  \item [$i)\Rightarrow ii)$]It suffices to prove that  $supp(\mu)\subset Z(P)$.\\
  Write  $\gamma_{ij}=\int\overline{z}^iz^jd\mu$, for  $i,j\geq0$.
   Since, for every $i\geq 0$   and  $r \le n$, we have
  $$\gamma_{i,n+1} -\sum\limits_{0\leq l+k\leq r} a_{lk} \gamma_{l+i,n+k-r}= 0,$$
  we get  $$\int\overline{z}^iz^{n+1} -\sum\limits_{0\leq l+k\leq r} a_{lk}\overline{z}^{l+i}z^{n+k-r} d\mu= 0.$$

  Hence  $$\int\overline{z}^iz^{n-r}P(z,\overline{z}) d\mu=0~~~~~~~~~~~~  \text{for every }  i\geq 0  \mbox{ and }  r \le n.$$
  Taking  an adequate combination, we manage to obtain   $$\int \overline{P}(z,\overline{z})P(z,\overline{z}) d\mu= \int \mid P(z,\overline{z})\mid^2 d\mu= 0.$$
  It follows that $P.\mu=0$, and thus $supp(\mu)\subset Z(P)$.

  \item [$ii)\Rightarrow i)$] Suppose that $\gamma_{ij}=\int\overline{z}^iz^jd\mu$, for all integers  $i, j$ such that $0\leq j\leq i\leq r$,
  and that $\text{supp}(\mu)\subset Z(P)$.
  Since $\gamma_{ij}=\overline{\gamma}_{ji}$, we also have $\gamma_{ij}=\int\overline{z}^iz^jd\mu$, for every  $i, j$ such that $0\leq i,j\leq r$.
Now
  \begin{equation*}
  \gamma_{i,r+1}=\sum\limits_{0\leq l+k\leq r} a_{lk}\gamma_{l+i,k}
  =\int\overline{z}^i(\sum\limits_{0\leq l+k\leq r} a_{lk}\overline{z}^lz^k)d\mu,
  \end{equation*}
and   since $\text{supp}(\mu)\subset Z(z^{r+1}-\sum\limits_{0\leq l+k\leq r} a_{lk}\overline{z}^lz^k)$,\\
  we get  $\gamma_{i,r+1}=\int\overline{z}^iz^{r+1}d\mu$.\\
  By induction we obtain $\gamma_{ij}=\int\overline{z}^iz^jd\mu$~~~~~~~~(for all $i,j\geq0$) and consequently, $\mu$ is a solution of the FCMP for $\{\gamma_{ij}\}_{i,j\geq0}$.
  \end{itemize}
\fin

\section{ Recursive Sequences of Fibonacci Type}
 In this section, we focus ourself on a particular case of RDIS, which will play a crucial role in the sequel.

 We shall refer to a RDIS with analytic characteristic polynomial as Recursive Sequences of Fibonacci Type (RSFT for short).
 In other words, a sequence $\gamma\equiv\{\gamma_{ij}\}_{i,j\geq0}$, with $\gamma_{00}>0$, is said to be RSFT, associated with the characteristic polynomial
 $P(x)=x^r-a_0x^{r-1}-\ldots-a_{r-2}x-a_{r-1}$, if it verifies the following relations:

 \begin{itemize}
\item [$i)$] For all $i,j\geq0$,
 \begin{equation}\label{2.1}
 \gamma_{ji}=\overline{\gamma}_{ij}.
 \end{equation}
\item [$ii)$] For all $i$ and $n$ such that $0\leq i\leq r-1\leq n$,
 \begin{equation}\label{2.2}
 \gamma_{i,n+1}=a_0\gamma_{i,n}+a_{1}\gamma_{i,n-1}+\ldots+a_{r-1}\gamma_{i,n-r+1}.
 \end{equation}
\end{itemize}

 We denote by $\mathcal{A}_\gamma$ the family of analytic characteristic, monic, polynomial associated with $\gamma$.

\begin{remarks}
 \begin{itemize}
  \item $\mathcal{A}_\gamma \subseteq \mathcal{P}_\gamma$.
  \item A sequence $\gamma\equiv\{\gamma_{ij}\}_{i,j\geq0}$ is a RSFT if and only if $\mathcal{A}_\gamma \neq \emptyset$.
 \end{itemize}
\end{remarks}

 \begin{prop}\label{p.2}
 Let $\gamma\equiv\{\gamma_{ij}\}_{i,j\geq0}$ be a RSFT, satisfying the relations \eqref{2.1} and \eqref{2.2}. Then $\mathcal{A}_\gamma$ is a principal ideal of $\mathbb{C}[X]$.
\end{prop}

{\it Proof. }
 It's obvious that $\mathcal{A}_\gamma$ is an ideal of $\mathbb{C}[X]$. It suffices to show that there exists a unique analytic characteristic polynomial,
 $P_\gamma \in \mathcal{A}_\gamma$, with minimal degree and that every characteristic polynomial is a multiple of $P_\gamma$.
  Since, for every $j\geq0$, the polynomial $P(x)=x^r-a_0x^{r-1}-\ldots-a_{r-2}x-a_{r-1}= \prod\limits_{k=1}^n (x-\lambda_i)^{d_i}$ is a characteristic polynomial associated with the singly indexed Fibonacci sequences $\gamma_j=\{\gamma_{ij}\}_{i\geq0}$, then (by \cite[Proposition 3.1]{crz}) there exists a unique characteristic  polynomial of minimal degree $P_{\gamma_j}= \prod\limits_{i=0}^n (x-\lambda_i)^{\alpha_{ij}}$, that divides $P(x)$, associated with $\gamma_j$. Since  $0\leq \alpha_{ij} \leq d_i$, for every $j\geq 0$, then $P_\gamma =\bigwedge\limits_{j\geq0} P_{\gamma_j}=  \prod\limits_{i=0}^n (x-\lambda_i)^{\alpha_{i}}$,  where $\alpha_i= \max_{j\geq0} \alpha_{ij}$, is the smallest common multiple of $P_{\gamma_j}$ divides $P$ and  provides a positive answer to the proposition.
  \fin

  Let $\gamma\equiv\{\gamma_{ij}\}_{i,j\geq0}$ be a RSFT, we call $P_\gamma$ the minimal polynomial associated with $\gamma$. Below, we  associate
  with  every RSFT its minimal  polynomial.

  \begin{cor}
  Let $\gamma\equiv\{\gamma_{ij}\}_{i,j\geq0}$ be a RSFT, associated with the minimal analytic characteristic polynomial $P_\gamma \in \mathbb{C}[z]$. If $M(\gamma) \geq 0$, then $P_\gamma$ has distinct roots.
  \end{cor}
  {\it Proof.} Setting  $P_\gamma(z) = \prod\limits_{i=1}^r (z -\lambda_i)^{n_i}$ for the characteristic polynomial of $(\gamma)$, we get  $M(\gamma) \prod\limits_{i=1}^r (z -\lambda_i)^{n_i} = 0$ (by Lemma \ref{l.0}).  It follows that  $M(\gamma) \prod\limits_{i=1}^r (z -\lambda_i)^{m} = 0$, where
  $m = \max\limits_{i=1}^r n^i$, and then  from Proposition \ref{l.4}, we have  $M(\gamma) \prod\limits_{i=1}^r (z -\lambda_i) = 0$. Therefore, again  by Lemma \ref{l.0}, the polynomial $\prod\limits_{i=1}^r (z -\lambda_i)$ is a characteristic polynomial of $\gamma$ and divides $P_\gamma$. Since $P_\gamma$ is minimal, then
  $P_\gamma = \prod\limits_{i=1}^r (z -\lambda_i)$.
  \fin

  We next solve the complex moment problem for a RSFT. Consider the following quadratic forms
  \begin{equation*}\begin{aligned}
  &\mathbb{C}[z, \overline{z}]&\rightarrow& \mathbb{C}\\
  \varphi^\gamma: &\sum\limits_{0\leq i+j\leq n} a_{ij}\overline{Z}^iZ^j &\rightarrow& \sum\limits_{0\leq i+j,h+k\leq n} a_{ij}a_{hk}\gamma_{i+k,j+h}
  \end{aligned}
  \end{equation*}
 and $$\varphi_n^\gamma\equiv \varphi^\gamma_{\mid\mathbb{C}_n[\overline{Z},Z]}.$$
  Let $M(\gamma)$ and $M(n)(\gamma)$ be the  matrices associated with $\varphi^\gamma$ and $\varphi_n^\gamma$, respectively. We denote  $M(n)(\gamma)\in M_m(\mathbb{C})$,  where $m=m(n)=\frac{(n+1)(n+2)}{2}$. Let also $\{e_{ij}\}_{0\leq i+j\leq n}$ be the canonical basis of ${\mathbb C}^m$, that is,  $e_{ij}$ is
  the vector with $1$ in the $\overline{Z}^iZ^j$ entry and $0$ all other positions.

  The next proposition establishes a link between the  positivity of $\varphi^\gamma$ and that of $\varphi_n^\gamma$.
  \begin{prop}\label{p.3}
   Let $\gamma$ be a RSFT  such that  $\deg P_\gamma=r$ and let  $n \geq 2r-2$. Then  $\varphi_n^\gamma$ is positive semi-definite
   if and only if $\varphi^\gamma$ is positive semi-definite. Moreover $\text{rank }\varphi_n^\gamma =\text{rank }\varphi^\gamma$.
 \end{prop}
 {\it Proof. }
We only need to show the direct implication. To this aim,  we construct a matrix $W\in M_{m,n+2}({\mathbb C})$ such that the   successive rows  are  defined by:

\centerline{$\overline{Z}^iZ^{2r-1-i}=\sum_{j=0}^{r-1} a_j e_{(i, 2r-2-i-j)},~~~~~~~\text{ for all}~~0\leq i\leq r-1,$}

 and

\centerline{$\overline{Z}^i Z^{2r-1-i}=\sum_{j=0}^{r-1} \overline{a_j} e_{(i-j-1, 2r-1-i)},~~~~~~~\text{ for all}~~r\leq i\leq 2r-1.$}

 Clearly  $M(n+1)(\gamma)$ has the form   $\begin{pmatrix} M(n)(\gamma)&B\\ B^*&C\end{pmatrix} $, with  $B=M(n)(\gamma)W$ and $C=B^*W$. Since
 $M(n)(\gamma)\geq0$ then, by the Smul'jan's theorem (see \cite{Smu}), $M(n+1)(\gamma)\geq0$ and $rank~M(n)(\gamma)=rank~M(n+1)(\gamma)$. By induction
  we obtain  $M(\gamma)\equiv M(\infty) \geq0$ and thus $\varphi^\gamma$ is positive semi-definite, as desired. \fin

 We are able now to  give a necessary and sufficient condition for a RSFT to be a moment sequence.
\begin{theorem}\label{t.3}
 Let $\gamma\equiv\{\gamma_{ij}\}_{i,j\geq0}$ be a RSFT, and $P_\gamma$  of degree $r$ be its minimal analytic characteristic polynomial. Then $\gamma$ admits a representing measure $\mu$ if and only if $\varphi^\gamma_{2r-2}$ is positive semi-definite. Moreover $$supp(\mu)= Z(P_\gamma).$$
\end{theorem}
{\it Proof.}
 If $\gamma$ admits a representing measure $\mu$; then, for any $p \in \mathbb{C}[z,\overline{z}]$, $p^T M(\gamma) p = \int \mid p\mid^2 d\mu \geq 0$. Thus it follows that $\varphi^\gamma$ is positive semidefinite, and hence   $\varphi^\gamma_{2r-2}$ is positive semi-definite. Conversely, if $M(2r -2)(\gamma)\geq0$, then  using  Proposition \ref{p.3}, we get $M(\gamma) \equiv M(\infty)$ is a positive semidefinite. Now let
 $P_\gamma (z) = \prod\limits_{i=1}^{r} (z-\lambda_i)$ and let $L_{\lambda_j}= \tiny{\prod\limits_{\begin{matrix} 1\leq i\leq r \\ i\neq j \end{matrix}}} \frac{z -\lambda_i}{\lambda_j -\lambda_i} \in \mathbb{C}[z,\overline{z}]$ ($i=1 \cdots ,r $)
  be the interpolation polynomials at the points of $Z(P_\gamma)$. Since the polynomials $Q(z, \overline{z}) = z^n -\sum_{i=1}^r \lambda_i^n L_{\lambda_i}$ and
 $H(z, \overline{z}) = \overline{z}^m -\sum_{j=1}^r \overline{\lambda}_j^m \overline{L_{\lambda_j}}$  vanish at all points of $Z(P_\gamma)$,  we derive from  Lemma \ref{p.11} that,  $Q(z, \overline{z})$ and $H(z, \overline{z})$ are characteristic polynomials of $\gamma$. Hence Lemma \ref{l.0} ensures that
 $M(\gamma) Q(z, \overline{z}) = M(\gamma) H(z, \overline{z})=0$, thus $\gamma_{mn} $ ($m,n \in \mathbb{Z}_+$) can be expressed as follows:
 \begin{equation*}\begin{aligned}
 \gamma_{mn} &= (z^m)^T M(\gamma) z^n\\
  &= 1^T M(\gamma) (z^n \overline{z}^m), &\text{by applying \eqref{p3.1},}\\
  &= 1^T M(\gamma) ( \sum_{i=1}^r \lambda_i^n L_{\lambda_i} \overline{z}^m), &\text{ from Remark  \ref{rmk}-i) and Lemma \ref{l.0},} \\
  &= 1^T M(\gamma) ( \sum_{i=1}^r \lambda_i^n L_{\lambda_i} \sum_{j=1}^r \overline{\lambda}_j^m \overline{L_{\lambda_j}})\\
  &= 1^T M(\gamma) (\sum_{i,j=1}^r \lambda_i^n \overline{\lambda}_j^m L_{\lambda_i} \overline{L_{\lambda_j}})\\
  &= \sum_{i,j=1}^r \lambda_i^n \overline{\lambda}_j^m 1^T M(\gamma) (L_{\lambda_i} \overline{L_{\lambda_j}}).
 \end{aligned}
 \end{equation*}
 If $i\neq j$, then $L_{\lambda_i} L_{\lambda_j}$ vanishing at all points of $Z(P_\gamma)$, hence (in this case and from Lemma \ref{p.11} ) $L_{\lambda_i} L_{\lambda_j}$ is a characteristic polynomial of $\gamma$, that is, $M(\gamma) (L_{\lambda_i} \overline{L_{\lambda_j}})=0$. Thus
$$\gamma_{mn} = \sum_{i=1}^r \overline{\lambda}_i^m \lambda_i^n 1^T M(\gamma) (L_{\lambda_i} \overline{L_{\lambda_i}})
 =\sum_{i=1}^r \lambda_i^n \overline{\lambda}_i^m L_{\lambda_i}^T M(\gamma) L_{\lambda_i}.$$
 Since $M(\gamma) \geq 0$,
 we get  $c_i = L_{\lambda_i}^T M(\gamma) L_{\lambda_i} \geq 0$  for $i = 0, \ldots, r$. Therefore the measure $\mu = \sum_{i=1}^r c_i \delta_{\lambda_i}$ gives  a positive answer to the problem \eqref{i.0} associated with $\gamma$.\\
 It remains to show that $supp(\mu)= Z(P_\gamma)$, that is, $c_i \neq 0$ for all $i=0, \ldots, r$. Indeed, if $c_i =0$ for some $i\in \{1, \ldots, r\}$, then $L_{\lambda_i}^T M(\gamma) L_{\lambda_i} =0$. Since $M(\gamma) \geq 0$, then $M(\gamma) L_{\lambda_i} =0$  and from Lemma \ref{l.0},  $L_{\lambda_i}$ is a
 analytic characteristic polynomial of $\gamma$, and this is a contradiction, because  the degree of the polynomial $L_{\lambda_i}$ is less strictly than the degree of the minimal polynomial $P_\gamma$.

  \section{ Solving the complex moment problem for RDIS}
 Let $\gamma=\{\gamma_{ij}\}_{i,j\geq0}$ be a double  indexed recursive \textit{moment} sequence, associated with the characteristic polynomial
 $P(z,\overline{z})= z^{r}-\sum\limits_{0\leq l+k\leq r-1} a_{lk}\overline{z}^lz^{k}$. Proposition \ref{p.1} ensures the existence of
 a representing measure $\mu$ associated with $\gamma$ such that  $supp(\mu)\subseteq Z(P):=\{\lambda_0,\ldots,\lambda_n\}$. It follows that
 $Q(z)=\prod\limits_{i=1}^n (z-z_i)= z^n -a_1z^{n-1} -\dots -a_n$ is a characteristic polynomial associated with $\gamma$, and since for every  $i\geq0$ and $n\geq r$, we have
 $$\gamma_{i,n} -a_1\gamma_{i,n-1} -a_{2}\gamma_{i,n-2} -\ldots -a_{r}\gamma_{i,n-r}= \int\overline{z}^i z^{n-r}Q(z) d\mu= 0.$$
 We conclude that a double indexed recursive \textit{moment} sequence is a RSFT.

  In order to obtain necessary and sufficient condition for a RDIS  to be moment sequence, we need to find the smallest $n$ which verifies the following equivalence $$M(n)\geq 0\Longleftrightarrow M(\gamma) \equiv M(\infty)\geq 0.$$ Since $\gamma$ admits a characteristic polynomial of the form $z^r - P_{r-1}(z,\overline{z})$ (with $\deg P_{r-1}\leq r-1$)  the best known result in this direction is the one  given by Curto-Fialkow \cite[Theorem 3.1]{3} which guarantees the equivalence for every  $n$  satisfying  the inequality $r\leq [\frac{n}{2}]+1$.  In the next theorem we involve the characteristic polynomials $P(z,\overline{z})$ and $Q(z)$ (associated with $\gamma$) to refine this result.

 \begin{theorem}\label{w.t.1}
 Let $\gamma\equiv\{\gamma_{ij}\}_{i,j\geq0}$ be a RDIS associated with the characteristic polynomials
 $P(z,\overline{z})= z^{r}-\sum\limits_{0\leq l+k\leq r-1} a_{lk}\overline{z}^l z^k$ and
 $Q(z)=\prod\limits_{i=1}^n (z-\lambda_i)$, where $Z(P)=\{\lambda_1, \ldots,\lambda_n\}$, then:
\begin{itemize}
\item [$i)$] There exists a characteristic polynomial, $h(z,\overline{z})\in \mathbb{C}[z,\overline{z}]$, associated with $\gamma$ such that
 $Q(z)=h(z,\overline{z})+f_1(z,\overline{z})P(z,\overline{z})+f_2(z,\overline{z})\overline{P}(z,\overline{z})$, where $f_1,f_2 \in \mathbb{C}[z,\overline{z}]$,
 with $d_z (h)<r$ and $d_{\overline{z}} (h)<r$.

\item [$ii)$]
 Let $A_h$ be the set of  monomials  $\overline{z}^iz^j   $  in $ h$  such that $ i+j=d_h $, and denote  $c_1 :=\max\{k\mid\overline{z}^lz^k\in A_h\}$ and $c'_1 :=\max\{l\mid\overline{z}^lz^k\in A_h\}$. We also define  $c_2 :=\max\{k\mid\overline{z}^lz^k\in A_h\setminus \{ \overline{z}^{d_h-c_1}z^{c_1}\}\}$
 and
 $c'_2 :=\max\{l\mid\overline{z}^lz^k\in A_h\setminus \{ \overline{z}^{c'_1}z^{d_h-c'-1}\}\}$, if card $A_h\geq 2$.    We put for convenience  $c_2=c'_2=-\infty$, in the case where  $card(A_h)=1$.
 We finally  denote $c= sup(c_1,c'_1)$,   $\alpha_{c_1} = inf(r-c_1, c_1-c_2)$ and  $\alpha_{c_1'} = inf(r-c'_1, c'_1-c'_2)$. Then
 $$  M(\infty)\geq 0    \Longleftrightarrow   M(2r-2-\alpha_c) \geq 0.$$
\end{itemize}
\end{theorem}

{\it Proof.}
 $i)$ It obvious to show that there exists $h(z,\overline{z})$ such that  $Q(z)=h(z,\overline{z})+f_1(z,\overline{z})P(z,\overline{z})+f_2(z,\overline{z})\overline{P}(z,\overline{z})$, for some $f_1,f_2 \in \mathbb{C}[z,\overline{z}]$, with $d_z (h)<r$ and $d_{\overline{z}} (h)<r$.
Since $P,\overline{P} \text{ and } Q$ are characteristic polynomials then $h(z,\overline{z})$ is also a characteristic polynomia associated with $(\gamma)$.

$ii)$ Since the two conditions are symmetric, we only need to prove the case  $c_1\geq c_1'$.
 Recall first that for $m\geq0$, $M_m(\mathbb{C})$ denotes the algebra of $m\times m$ complex matrices and for $n\geq0$, let $m\equiv m(n) := (n+1)(n+2)/2$ and $M(n)(\gamma)\in M_m(\mathbb{C})$, as in the introduction. We define a basis $\{e_{ij}\}_{0\leq i+j\leq n}$ of ${\mathbb C}^m$ as follows: $e_{ij}$ is the vector with $1$ in the $\overline{Z}^iZ^j$ entry and $0$ all other positions.

 The main idea is to write, from the characteristic polynomials associated with $\gamma$, monomials of order $2r-1-\alpha_c+e-1$ ($e\in \mathbb{N}$) as linear combination of monomials of order strictly less than $2r-1-\alpha_c+e$ modulo $P,\overline{P},h$ and $\overline{h}$, that is,
 $$\overline{z}^i z^{2r-1-\alpha_c+e-i}=H_{(i,2r-1-\alpha_c+e-i)}(z,\overline{z}) \text{ mod }\{P,\overline{P},h,\overline{h}\},$$
 $\text{ with } 0\leq i\leq 2r-1-\alpha_c+e \text{ and } \deg H_{(i, 2r-1-\alpha_c+e-i)}(z,\overline{z})\leq 2r-2-\alpha_c.$

 We construct for every $e \in \mathbb{Z}_+$ a matrix $W_e \in M_{m(2r-2-\alpha_c+e), 2r-1-\alpha_c+e}(\mathbb{C})$ such that the coefficients of the column $\overline{Z}^i Z^{2r-1-\alpha_c+e-i}$ , $0\leq i \leq  2r-1-\alpha_c+e$, are that of the polynomial $H_{(i,2r-1-\alpha_c+e-i)}$ (considering the lexicographic order cited in the  introduction).

 Since $P,\overline{P},h$ and $\overline{h}$ are characteristic polynomials associated with $\gamma$,   the above discussion leads, in view of Lemma \ref{l.0}, to the following equality:
 $$M(2r-1-\alpha_c+e)(\gamma)=  \begin{pmatrix} M(2r-2-\alpha_c+e)(\gamma)&B\\ B^*&C\end{pmatrix},$$
 such that $B=M(2r-2-\alpha_c+e)(\gamma)W_e$ and $C=B^*W_e$, for all $e\in \mathbb{N}$.
 According to Smul'jan's theorem we have $$M(2r-2-\alpha_c+e)(\gamma)\geq0 \Longleftrightarrow M(2r-1-\alpha_c+e)(\gamma)\geq0,$$  and hence it follows by induction that
 $$M(2r-2-\alpha_c)(\gamma)\geq0 \Longleftrightarrow M(\infty)(\gamma)\geq0.$$
  We distinguish 3  cases,

  \begin{itemize}
   \item  $r\leq i\leq 2r-1-\alpha_c+e$.
 Since $\overline{z}^r=\overline{P} +\sum\limits_{0\leq l+k\leq r-1} \overline{a}_{lk}\overline{z}^k z^l$, we obtain
 \begin{equation*}
 \overline{z}^i z^{2r-1-\alpha_c+e-i} = \overline{z}^{i-r}z^{2r-1-\alpha_c+e-i}\overline{P}+ \sum\limits_{0\leq l+k\leq r-1} \overline{a}_{lk}\overline{z}^{k+i-r} z^{l+2r-1-\alpha_c+e-i}.
 \end{equation*}
 Hence
 \begin{equation*}
 \overline{Z}^i Z^{2r-1-\alpha_c+e-i}=  \sum\limits_{0\leq l+k\leq r-1} \overline{a}_{lk} e_{(k+i-r, l+2r-1-\alpha_c+e-i)}.
 \end{equation*}

   \item $0\leq i\leq r-1-\alpha_c.$
   As  $z^r= P +\sum\limits_{0\leq l+k\leq r-1} a_{lk}\overline{z}^l z^k$, we get
 \begin{equation*}
 \overline{z}^i z^{2r-1-\alpha_c+e-i}= \overline{z}^{i} z^{r-1-\alpha_c+e-i}P + \sum\limits_{0\leq l+k\leq r-1} a_{lk}\overline{z}^{l+i} z^{k+r-1-\alpha_c+e-i},
 \end{equation*}
 and thus,
 \begin{equation*}
 \overline{Z}^i Z^{2r-1-\alpha_c+e-i} =\sum\limits_{0\leq l+k\leq r-1} a_{lk} e_{(l+i, k+r-1-\alpha_c+e-i)}.
 \end{equation*}

   \item $r-\alpha_c \leq i\leq r-1$. This third case requires more work. We will distinguish two sub cases.\\

   $a)$  Card $A=1$. (In this case we have $\alpha_c = r-c_1$).
    Let  $$h(z, \overline{z}) =\overline{z}^{d_h -c_1} z^{c_1} -\sum\limits_{0\leq l+k\leq d_h -1} a_{lk}\overline{z}^lz^k,$$
    we have
   $$\overline{z}^{d_h -c_1} z^{c_1}= h(z, \overline{z}) +\sum\limits_{0\leq l+k\leq d_h -1} a_{lk}\overline{z}^l z^k.$$
    and  since $ d_h \leq c_1 +c'_1 \leq 2 c_1,$ we get
   \begin{equation*}\begin{array}{lll}
 \overline{z}^{r+c_1+e-1-i} z^i &=& \overline{z}^{r+ 2c_1+e-1-i-d_h} z^{i-c_1}(\overline{z}^{d_h -c_1} z^{c_1}),\\
  &=& \overline{z}^{r+ 2c_1+e-1-i-d_h} z^{i-c_1} h(z, \overline{z}) \\
  & &+\sum\limits_{0\leq l+k\leq d_h -1} a_{lk}\overline{z}^{l+r+ 2c_1+e-1-i-d_h} z^{k+i-c_1}.
  \end{array}
 \end{equation*}

 Hence,
 \begin{equation*}
 \overline{Z}^{r+c_1+e-1-i} Z^i=\sum\limits_{0\leq l+k\leq d_h -1} a_{lk}e_{(l+r+ 2c_1+e-1-i-d_h, k+i-c_1)},
 \end{equation*}
and it follows that,
 \begin{equation*}
 \overline{Z}^{2r-1-\alpha_c+e-i} Z^i=\sum\limits_{0\leq l+k\leq d_h -1} a_{lk}e_{(l+3r-2\alpha_c+e-1-i-d_h, k+i-r+\alpha_c)}.
 \end{equation*}

 $b)$ Card $A\geq2$. Let
 $$\overline{z}^{d_h-c_1} z^{c_1} =h(z, \overline{z}) +\sum\limits_{l+k=d_h} \alpha_{lk} \overline{z}^lz^k +\sum\limits_{0\leq l+k\leq d_h-1}a_{lk}\overline{z}^lz^k,$$
 we  deduce that,
 $$\begin{array}{ll}
 \overline{z}^{c_2 -c_1 +r} z^{c_1} &= (\overline{z}^{r-d_h+c_2}) (\overline{z}^{d_h-c_1}z^{c_1})\\
 &=\overline{z}^{r-d_h+c_2}h+\sum\limits_{l+k=d_h} \alpha_{lk} \overline{z}^{l+r-d_h+c_2}z^k\\
  & +\sum\limits_{0\leq l+k\leq d_h-1}a_{lk}\overline{z}^{l+r-d_h+c_2}z^k.
   \end{array}$$
 Since in the above equality,  all monomials in the sum $\sum\limits_{l+k=d_h} \alpha_{lk} \overline{z}^{r+l-d_h+c_2}z^k$ satisfy
 $d_{\overline{z}} (\overline{z}^{l+r-d_h+c_2}z^k) \geq r$ (since $k\leq c_2$ and $l+k=d_h$, then $l -d_h +c_2\geq 0$),  we get
 $$\sum\limits_{l+k=d_h} \alpha_{lk} \overline{z}^{l+r-d_h+c_2}z^k =  \overline{P}\sum\limits_{l+k=d_h} \overline{z}^{l-d_h+c_2}z^k - \sum\limits_{0\leq l'+k'\leq r+c_2-1} a_{l'k'} \overline{z}^{l'}z^{k'},$$
 where $\{a_{l'k'}\}_{0\leq l'+k'\leq r+c_2-1}$ are complex numbers. Thus, there exists $ \alpha_{l''k''} \in \mathbb{C}$ such that
 \begin{equation}\label{pr.1}
 \overline{z}^{ c_2-c_1+r}z^{c_1} = \sum\limits_{0\leq l''+k''\leq r+c_2-1} \alpha_{l''k''}\overline{z}^{l''}z^{k''} mod(\overline{P},h).
 \end{equation}
Here again we discuss two situations,\\
 *) If $r-c_1 \leq c_1-c_2$, (that is, $\alpha_c = r-c_1$).\\
   Then $0\leq c_1 -r +c_1 -c_2 \leq i -r +c_1 -c_2$ and $0\leq c_2-c_1+r$ (since $c_2\leq c_1\leq i\leq r$), hence \eqref{pr.1} yields
 \begin{equation*}\begin{aligned}
 \overline{z}^i z^{r+c_1+e-1-i} &= (\overline{z}^{i-c_2+c_1-r}z^{r+e-1-i}) (\overline{z}^{c_2-c_1+r}z^{c_1})\\
 &=\sum\limits_{0\leq l''+k''\leq r+c_2-1} \alpha_{l''k''}\overline{z}^{l''+i-c_2+c_1-r}z^{k''+r+e-1-i} mod(\overline{P},h).
 \end{aligned}
 \end{equation*}
 Then  $$\overline{Z}^i Z^{r+c_1+e-1-i} = \sum\limits_{0\leq l''+k''\leq r+c_2-1} \alpha_{l''k''} e_{(l''+i-c_2+c_1-r, k''+r+e-1-i)},$$
 that is, $$\overline{Z}^i Z^{2r-1-\alpha_c+e-i} = \sum\limits_{0\leq l''+k''\leq r+c_2-1} \alpha_{l''k''} e_{(l''+i-c_2+c_1-r, k''+r+e-1-i)}.$$
 **) If $r-c_1> c_1-c_2$, (that is, $\alpha_c = c_1-c_2$).\\
  By  \eqref{pr.1}, we get
 \begin{equation*}\begin{aligned}
 \overline{z}^{2r-c_1+c_2+e-1-i} z^i &= (\overline{z}^{r+e-1-i} z^{i-c_1}) (\overline{z}^{r-c_1+c_2} z^{c_1})\\
 &=\sum\limits_{0\leq l''+k''\leq r+c_2-1} \alpha_{l''k''}\overline{z}^{l''+r+e-1-i}z^{k''+i-c_1} mod(\overline{P},h),
 \end{aligned}
 \end{equation*}
 and then

 $$\overline{z}^i z^{2r-1-\alpha_c+e-i} = \sum\limits_{0\leq l''+k''\leq r+c_2-1} \overline{\alpha}_{l''k''} \overline{z}^{k''+i-c_1} z^{l''+r+e-1-i} mod(P,\overline{h}).$$
 This finally gives $$\overline{Z}^i Z^{2r-1-\alpha_c+e-i}= \sum\limits_{0\leq l''+k''\leq r+c_2-1} \overline{\alpha}_{l''k''} e_{(k''+i-c_1, l''+r+e-1-i)},$$ as required.
 \end{itemize}
 \fin

 \begin{cor}\label{c.t.1}
 Let $\gamma\equiv\{\gamma_{ij}\}_{i,j\geq0}$ be a RDIS associated with the characteristic polynomial
 $P(z,\overline{z})= z^{r}-\sum\limits_{0\leq l+k\leq r-1} a_{lk}\overline{z}^l z^k$, then:
 $$  M(\infty)(\gamma)\geq 0    \Longleftrightarrow   M(2r-2)(\gamma) \geq 0.$$
\end{cor}

 {\it Proof.}
 It suffices to rewrite the first and second cases, in the above theorem's proof, with $\alpha_c =0$.
 \fin

 In the next theorem and in the sequel, we set  $\xi\equiv\xi_\gamma = 2r-2-\alpha_c$.

\begin{theorem}\label{t.2}
 Let $\gamma\equiv\{\gamma_{ij}\}_{i,j\geq0}$ be a RDIS associated with the characteristic polynomial
 $P(z,\overline{z})= z^{r}-\sum\limits_{0\leq l+k\leq r-1} a_{lk}\overline{z}^lz^{k}$ and let $Z(P)= \{\lambda_1,\lambda_2,\ldots,\lambda_n\}$. The following are equivalent:
\begin{itemize}
\item $\gamma$ is a moment sequence;
\item  $Q(z)=\Pi_{i=1}^n(z-\lambda_i) \in \mathcal{P}_\gamma$ and $M(\xi)\geq0$.
\end{itemize}
\end{theorem}

{\it Proof.}
  Suppose that $\gamma$ admits a representing measure $\mu$.\\
Then, from Proposition \ref{p.1}, $\text{supp}(\mu)\subset Z(P)= \{\lambda_1,\lambda_2,\ldots,\lambda_n\}$.\\
Setting  $Q(z)=\Pi_{i=1}^n(z-\lambda_i)=z^n-\alpha_1z^{n-1}-\ldots-\alpha_n$, we have
\begin{equation*}\begin{aligned}
\gamma_{im}&=\int \overline{z}^iz^m d\mu,~~~~~~~~~~~~~\text{ where }~~0\leq i \text{ and } n\leq m,&\\
&=\int\overline{z}^i(\alpha_1z^{m-1}+\alpha_2z^{m-2}+\ldots+\alpha_n z^{m-n})d\mu&\\
&=\alpha_1\gamma_{i,m-1}+\alpha_2\gamma_{i,m-2}+\ldots+\alpha_n \gamma_{i,m-n}.&
\end{aligned}\end{equation*}
Hence $Q(z)$ is a characteristic polynomial associated with $\gamma$. The condition $M(\xi)\geq0$ is obvious. Conversely, since  $M(\gamma)(\xi)\geq0$ then, from  Theorem \ref{w.t.1},
   $M(\gamma)(\infty)\geq0$ and thus, from Theorem \ref{t.3}, the sequence $\gamma$ owns a representing measure.
\fin

 \begin{cor}\label{c.t.2}
 Let $\gamma\equiv\{\gamma_{ij}\}_{i,j\geq0}$ be a RDIS associated with the characteristic polynomial
 $P(z,\overline{z})= z^{r}-\sum\limits_{0\leq l+k\leq r-1} a_{lk}\overline{z}^lz^{k}$, then $\gamma$ is a moment sequence if and only if
 $M(2r-2)(\gamma) \geq0$.
 \end{cor}
 {\it Proof.}
 We only need to show the converse implication. As noted in the introduction the polynomial $P(z,\overline{z})$ has finite number of roots, say $Z(P) := \{ \lambda_1, \lambda_2, \ldots, \lambda_n\}$. Since $M(2r-2)(\gamma) \geq0$,  Corollary \ref{c.t.1} yields $M(\infty)(\gamma) \geq 0$. Let $Q(z)=\Pi_{i=1}^n (z-\lambda_i)$, applying Proposition \ref{p.11},  we obtain $Q(z,\overline{z}) \in \mathcal{P}_\gamma$. Therefore, Theorem \ref{t.2} implies that $\gamma$ is a moment sequence.
 \fin

\section{ The case of cubic moment problem }
\subsection{The TCMP with cubic relation of the form $z^3+az+b\overline{z}=0$}

Whenever we have the  zeros of a characteristic  polynomial associated with a RDIS,  Theorem \ref{t.2} allows us to give a concrete, computable, necessary and
 sufficient conditions for the existence of a representing measure associated with this sequence. In this section we apply the above mentioned theorem to solve
 the complex moment problem for a RDIS, associated with a harmonic characteristic polynomial of the form $z^3+az+b\overline{z}$, where $a,b \in \mathbb{R}$.
 First, we start by giving the number of zeros of the harmonic polynomial $P(z,\overline{z})=z^3+az+b\overline{z}$. We solve the equation  $z^3+az+b\overline{z}=0$, for completeness.
 Writing $z=x+iy$, we get   $(x+iy)^3+a(x+iy)+b(x-iy)=0$, and
 then $x^3-3xy^2+(a+b)x-i(y^3-(3x^2+a+d)y)=0$. It follows that
 \begin{equation*}
 \begin{cases}
 x(x^2-3y^2+a+b)&=0,\\
 y(y^2-(3x^2+a-b))&=0.
 \end{cases}
\end{equation*}

\begin{itemize}
  \item If $y=0$; then $x=0$ or $x^2=-a-b$.
  \item If $x=0$; then $y=0$ or $y^2=a-b$.
  \item If $xy\neq 0$; then $x^2=\frac{-a+2b}{4}$ and $y^2=\frac{a+2b}{4}$.
\end{itemize}
We deduce  from the above cases that $z^3+az+b\overline{z}$ has at most  7 roots and it has exactly 7 roots  if and only if $b< |a|<2b$. \\

Let $\lambda_1,\lambda_2,\ldots,\lambda_7$ be the roots of the polynomial $P(z,\overline{z})=z^3+az+b\overline{z}$, with $b<-a<2b$. Direct computations lead to the expression
\begin{equation*}\begin{aligned}
Q(z)&=\Pi_{i=1}^7(z-\lambda_i)&\\
&=z^7+(2a+b)z^5+(a^2+b^2+ab)z^3+(b^3+ab^2)z&\\
&=(z^4+(a+b)z^2-b\overline{z}z+b^2)P(z,\overline{z}) -b^2 (\overline{z}z^2 -\overline{z}^2z -bz +b\overline{z})&\\
&=(z^4+(a+b)z^2-b\overline{z}z+b^2)P(z,\overline{z}) -b^2 h(z, \overline{z}).&
\end{aligned}\end{equation*}

 Now we are able to solve  the moment problem of this section, in the case $b<-a<2b$. Recall again that,  $1,Z,\overline{Z},\dots,Z^n,\ldots,\overline{Z}^n$  denote the successive columns of $M(\omega)(n)$.

\begin{theorem}\label{t.4}
 Let $\omega\equiv\{\gamma_{ij}\}_{0\leq i+j\leq 6}$, with $\overline{\gamma}_{ij}=\gamma_{ji}$ and $\gamma_{00}>0$, be a truncated complex sequence, and let
 $M(\omega)(3)$ be its associated moment matrix. If $M(\omega)(3)\geq 0$ and has cubic column relation of the form $Z^3=-aZ-b\overline{Z}$, with $a,b\in \mathbb{R}$ and $b<-a<2b$.
 Then the following statements are equivalent:
\begin{itemize}
  \item [$i)$] There exists a representing measure for $\omega$.

  \item [$ii)$] There exists a representing measure for the RDIS $\gamma\equiv\{\gamma_{ij}\}_{i,j\geq 0}$, whose initial conditions and characteristic
  polynomial are $\{\gamma_{ij}\}_{0\leq i\leq j\leq 2}$ and $P(z,\overline{z})=z^3+az+b\overline{z}$, respectively .

  \item [$iii)$]
    \begin{equation*}\begin{cases}
\Lambda_\omega(h)&=0,\\
\Lambda_\omega(zh)&=0.
    \end{cases}\end{equation*}

  \item [$iv)$]
    \begin{equation*}\begin{cases}
\mathfrak{Im}\gamma_{12}&=b\mathfrak{Im}\gamma_{01},\\
\gamma_{22}+2b\mathfrak{Re}\gamma_{20}+(a-b)\gamma_{11}&=0.
    \end{cases}\end{equation*}

  \item [$v)$] $h(z,\overline{z})=\overline{z}z^2 -\overline{z}^2z -bz +b\overline{z} \in \mathcal{P}_\gamma$, where $\gamma$ is the RDIS defined in $ii)$.
  \item [$vi)$]  $\overline{Z}Z^2 -\overline{Z}^2Z -bZ +b\overline{Z}=0$.
\end{itemize}
\end{theorem}

{\it proof.}
It is straightforward to see that $iii)\Leftrightarrow iv)$, $ii)\Rightarrow i)$ and $v)\Rightarrow vi) \Rightarrow iii)$.
Thus, it is enough to show $i)\Rightarrow ii)$, $iii)\Rightarrow v)$ and $ v)\Leftrightarrow ii)$.\\

  $i)\Longrightarrow ii)$.  Suppose that $\omega$ admits a representing measure $\mu$. Because the cubic column in $M(\omega)(3)$ verifies $Z^3= -aZ-b\overline{Z}$, we get  the following relations
 \begin{equation*}\begin{aligned}
 \int \overline{z}^3 P(z,\overline{z}) d\mu &=\int \overline{z}^3 (z^3+az+b\overline{z}) d\mu &= \gamma_{3,3} +a\gamma_{3,1} +b\gamma_{4,0} &=0,\\
 \int z P(z,\overline{z}) d\mu &=\int z (z^3+az+b\overline{z}) d\mu &= \gamma_{0,4} +a\gamma_{0,2} +b\gamma_{1,1} &=0,\\
 \int \overline{z} P(z,\overline{z}) d\mu &=\int \overline{z} (z^3+az+b\overline{z}) d\mu &= \gamma_{1,3} +a\gamma_{1,1} +b\gamma_{2,0} &=0,
 \end{aligned}\end{equation*}
 thus $\int \mid P\mid^2 d\mu = \int (\overline{z}^3+a\overline{z}+b z)P d\mu =0$ hence $supp(\mu)\subseteq Z(P)$. It follows, from Proposition \ref{p.1}, that $\mu$ is a representing measure for $\gamma$.\\
 $iii)\Rightarrow v)$. It suffices to show that $\Lambda_\gamma(\overline{z}^iz^jh(z,\overline{z}))=0$, for all $i,j\geq0$.\\
 Remark that, whenever $z^3 = -az -b\overline{z}$, we have:\\
 \begin{equation*}\begin{aligned}
&z^2h&=&(b-a)h,&\\
&\overline{z}zh&=&(a-b)h,&\\
&\overline{h}&=&-h.&
\end{aligned}\end{equation*}
Then
 \begin{equation*}\begin{array}{l}
\Lambda_\gamma(z^2h)                    =(b-a)\Lambda_\omega(h)                 = 0,\\
\Lambda_\gamma(\overline{z}zh)     =(a-b)\Lambda_\omega(h)                 = 0,\\
\Lambda_\gamma(\overline{h})        =-\Lambda_\omega(h)                         = 0,\\
\Lambda_\gamma(\overline{z}h)       =\overline{\Lambda_\omega(z\overline{h})}=-\overline{\Lambda_\omega(zh)}=0,\\
\Lambda_\gamma(\overline{z}z^2h) =(b-a)\Lambda_\omega(\overline{z}h)         =-(b-a)\overline{\Lambda_\omega(zh)}=0.
\end{array}\end{equation*}

 Since $\gamma$ is a RDIS associated with the characteristic polynomial $z^3+az+b\overline{z}$, then, for every $i,j \in \mathbb{Z}_+$, we have
$$\Lambda_\gamma(\overline{z}^iz^jh(z,\overline{z})) = \sum\limits_{0\leq l,k \leq 2} a_{lk} \Lambda_\omega(\overline{z}^l z^k h(z,\overline{z})) =0,$$ where $\{a_{lk}\}_{0\leq l,k \leq 2}$ are real numbers.

$v)\Rightarrow ii)$.  We know that $Q(z)=(z^4+(a+b)z^2-b\overline{z}z+b^2)P(z,\overline{z})-b^2h(z,\overline{z})$.  Since $P(z,\overline{z})$ and $h(z,\overline{z})$ are  characteristic polynomials, $Q(z)$ is also a characteristic polynomial associated with $\gamma$. As $M(\xi)(\gamma) \equiv M(3) \geq 0$ (observe that $\xi = 2\times 3-2+1 =3$) and $Q(z)\in \mathcal{P}_\gamma$, then Theorem \ref{t.2} yields that $\gamma$ admits a representing measure.

$ii)\Rightarrow v)$.   We have
$$ h(z,\overline{z})= \frac{1}{b^2 }(z^4+(a+b)z^2-b\overline{z}z+b^2)P(z,\overline{z}) -\frac{1}{b^2 } Q(z),$$
since (by Theorem \ref{t.2}) $Q(z)$ is a characteristic polynomial of $\gamma$, as well as the polynomial $P(z,\overline{z})$, then
$h(z,\overline{z})$ is a characteristic polynomial associated with $\gamma$ (by Remark \ref{rmk}-$i)$ ).
\fin

 Let us now suppose that $b<a<2b$, then
\begin{equation*}\begin{aligned}
Q(z)&=\Pi_{i=1}^7(z-\lambda_i)&\\
&=z^7+(2a-b)z^5+(a^2+b^2-ab)z^3+(ab^2-b^3)z&\\
&=(z^4+(a-b)z^2-b\overline{z}z+b^2)P(z,\overline{z})+b^2(z+\overline{z})(\overline{z}z-b).&
\end{aligned}\end{equation*}

Let $h(z,\overline{z})=(z+\overline{z})(\overline{z}z-b)$, we have
 \begin{equation*}\begin{aligned}
z^2h(\overline{z},z)&=(\overline{z}z+\overline{z}^2-b)P(z,\overline{z})-b\overline{P}(\overline{z},z)-(a+b)h(\overline{z},z),&\\
z\overline{z}h(\overline{z},z)&=(z^2-b)\overline{P}(\overline{z},z)+(\overline{z}^2-b)P(\overline{z},z)-(a+b)h(\overline{z},z),&\\
\overline{h}&=h.&
\end{aligned}\end{equation*}
 If $\Lambda_\omega(h) = \Lambda_\omega(zh) =0$, we get
 \begin{equation*}\begin{aligned}
&\Lambda_\gamma(\overline{h})&=& \Lambda_\omega(h) &=0,&\\
&\Lambda_\gamma(\overline{z}h)&=&\overline{\Lambda_\gamma(z\overline{h})}=\overline{\Lambda_\omega(zh)} &=0,&\\
&\Lambda_\gamma(z\overline{z}h)&=&-(a+b)\Lambda_\omega(h) &=0,&\\
&\Lambda_\gamma(z^2h)&=&-(a+b)\Lambda_\omega(h) &=0,&\\
&\Lambda_\gamma(\overline{z}z^2h)&=&-(a+b)\overline{\Lambda_\omega(zh)} &=0,&\\
&\Lambda_\gamma(\overline{z}^2zh)&=&-(a+b)\Lambda_\omega(zh) &=0.&\\
\end{aligned}\end{equation*}
 The above equalities and cubic relation allow us to prove that $\Lambda_\gamma(z^i\overline{z}^jh)=0$, for all $i,j\in\mathbb{Z}_+$; that is, $h(\overline{z},z) \in \mathcal{P}_\gamma$. From this observation  and similarly to the above theorem's proof, we are able to state the following theorem.

\begin{theorem}\label{t.5}
 Let $\omega\equiv\{\gamma_{ij}\}_{0\leq i+j\leq 6}$, with $\overline{\gamma}_{ij}=\gamma_{ji}$ and $\gamma_{00}>0$, be a given truncated complex sequence, and let
 $M(\omega)(3)$ be its associated moment matrix. If $M(\omega)(3)\geq 0$ and has cubic column relation of the form $Z^3=-aZ-b\overline{Z}$ , with $a,b\in \mathbb{R}$ and $b<a<2b$.
 Then the following statements are equivalent:
\begin{itemize}
  \item [$i)$] There exists a representing measure for $\omega$.

  \item [$ii)$] There exists a representing measure for the RDIS $\gamma\equiv\{\gamma_{ij}\}_{i,j\geq 0}$, whose initial conditions and characteristic
  polynomial are $\{\gamma_{ij}\}_{0\leq i\leq j\leq 2}$ and $P(z,\overline{z})=z^3+az+b\overline{z}$, respectively .

  \item [$iii)$]
    \begin{equation*}\begin{cases}
\Lambda_\omega(h)&=0,\\
\Lambda_\omega(zh)&=0.
    \end{cases}\end{equation*}

  \item [$iv)$]
    \begin{equation*}\begin{cases}
\mathfrak{Re}\gamma_{12}&=b\mathfrak{Re}\gamma_{01},\\
\gamma_{22}&=2b\mathfrak{Re}\gamma_{20}+(a+b)\gamma_{11}.
    \end{cases}\end{equation*}

  \item [$v)$] $h(z,\overline{z})=(z+\overline{z})(z\overline{z}-b)\in \mathcal{P}_\gamma$, where $\gamma$ is the RDIS defined in $ii)$.
  \item [$vi)$]  $(Z+\overline{Z})(Z\overline{Z}-b)=0$.
\end{itemize}
\end{theorem}

By the same technique we treat the other cases, we find the following results.
\small{
\begin{table}[htbp]
\begin{center}
\begin{tabular}{|c|c|c|c|c|c|}
  \cline{2-6}
  \multicolumn {1}{ c | }{$b>0$ } & N & Zeros & $h(z,\overline{z})$ & \multicolumn {2}{| c | }{Necessary and sufficient conditions } \\

  \hline

  $2b\leq a$ & 3 &
  ~$\begin{aligned} \pm& i\sqrt{a-b}\\ &0 \end{aligned}$~&
  ~$z+\overline{z}$~&
  ~$\begin{aligned} &M(2)\geq0 \\ &\Lambda_\gamma(h)=0 \\ &\Lambda_\gamma(zh)=0 \\ &\Lambda_\gamma(z^2h)=0 \end{aligned}$~&
  ~$\begin{aligned} &M(2)\geq0 \\ &\mathfrak{Im}\gamma_{01}=0 \\ &\gamma_{11}+\gamma_{02}=0 \\ &\gamma_{12}=a\gamma_{01}+b\gamma_{10} \end{aligned}$~\\

  \hline

  $b<a<2b$ & 7 &
  ~$\begin{aligned} \pm& i\sqrt{a-b};\\ \pm& \frac{\sqrt{-a+2b}}{2}\\ \pm& i\frac{\sqrt{a+2b}}{2};\\ &0. \end{aligned}$~&
  ~$(z+\overline{z})(z\overline{z}-u)$~&
  ~$\begin{aligned} &M(3)\geq0 \\ &\Lambda_\gamma(h)=0 \\ &\Lambda_\gamma(zh)=0 \end{aligned}$~&
  ~$\begin{aligned} &M(3)\geq0 \\ &\mathfrak{Re}\gamma_{12}=b\mathfrak{Re}\gamma_{01}\\ &\gamma_{22}=2b\mathfrak{Re}\gamma_{20}+(a+b)\gamma_{11} \end{aligned}$~\\

  \hline

  $-b\leq a\leq b$ & 1 &
  ~0~&
  ~ ~&
  ~ ~&
  ~$\begin{aligned} &\gamma_{00}>0 \text{ and } \gamma_{ij}=0 \\ &\text{for all } 0\leq i\leq j\leq 2 \end{aligned}$~\\

  \hline

  $b<-a<2b$ & 7 &
  ~$\begin{aligned} \pm& \sqrt{-a-b};\\ \pm& \frac{\sqrt{-a+2b}}{2}\\ \pm& i\frac{\sqrt{a+2b}}{2};\\ &0. \end{aligned}$~&
  ~$(z-\overline{z})(z\overline{z}-u)$~&
  ~$\begin{aligned} &M(3)\geq0 \\ &\Lambda_\gamma(h)=0 \\ &\Lambda_\gamma(zh)=0 \end{aligned}$~&
  ~$\begin{aligned} &M(3)\geq0 \\ &\mathfrak{Im}\gamma_{12}=b\mathfrak{Im}\gamma_{01}\\ &\gamma_{22}+2b\mathfrak{Re}\gamma_{20}=(b-a)\gamma_{11} \end{aligned}$~\\

  \hline

  $a\leq-2b$ & 3 &
  ~$\begin{aligned} \pm& \sqrt{-a-b}\\ &0 \end{aligned}$~&
  ~$z-\overline{z}$~~&
  ~$\begin{aligned} &M(2)\geq0 \\ &\Lambda_\gamma(h)=0 \\ &\Lambda_\gamma(zh)=0 \\ &\Lambda_\gamma(z^2h)=0 \end{aligned}$~&
  ~$\begin{aligned} &M(2)\geq0 \\ &\gamma_{01}=\gamma_{10}\\ &\gamma_{02}=\gamma_{11} \\ &a\gamma_{01}+b\gamma_{10}+\gamma_{12}=0 \end{aligned}$~\\

\hline
\end{tabular}
\end{center}
\caption{$b$ positive}
\label{tab:table1}
\end{table}
}
\begin{table}[htbp]
\begin{center}
\begin{tabular}{|c|c|c|c|c|c|}
  \cline{2-6}
  \multicolumn {1}{ c | }{$b<0$ } & N & Zeros & $h(z,\overline{z})$ & \multicolumn {2}{| c | }{Necessary and sufficient conditions } \\
  \hline
  $-b\leq a$ & 3 &
  ~~$\begin{aligned} \pm& i\sqrt{a-b}\\ &0 \end{aligned}$~~&
  ~~$z+\overline{z}$~~&
  ~~$\begin{aligned} &M(2)\geq0 \\ &\Lambda_\gamma(h)=0 \\ &\Lambda_\gamma(zh)=0 \\ &\Lambda_\gamma(z^2h)=0 \end{aligned}$~~&
  ~~$\begin{aligned} &M(2)\geq0 \\ &\mathfrak{Im}\gamma_{01}=0 \\ &\gamma_{11}+\gamma_{02}=0 \\ &\gamma_{12}=a\gamma_{01}+b\gamma_{10} \end{aligned}$~~\\

  \hline

  $|a|<-b$ & 5 &
  ~~$\begin{aligned} \pm& i\sqrt{a-b}\\ \pm& \sqrt{-a-b}\\ &0 \end{aligned}$~~&
  ~~$z^2\overline{z}+a\overline{z}+bz$~~&
  ~~$\begin{aligned} &M(3)\geq0 \\ &\Lambda_\gamma(h)=0 \\ &\Lambda_\gamma(zh)=0 \\ &\Lambda_\gamma(\overline{z}h)=0 \end{aligned}$~~&
  ~~$\begin{aligned} &M(3)\geq0 \\ &\gamma_{21}+a\gamma_{01}+b\gamma_{10}=0 \\ &\gamma_{20}=\gamma_{02} \\ &\gamma_{22}+a\gamma_{01}+b\gamma_{10}=0 \end{aligned}$~~\\

  \hline

  $a \leq b$ & 3 &
  ~~$\begin{aligned} \pm& \sqrt{-a-b}\\ &0 \end{aligned}$~~&
  ~~$z-\overline{z}$~~&
  ~~$\begin{aligned} &M(2)\geq0 \\ &\Lambda_\gamma(h)=0 \\ &\Lambda_\gamma(zh)=0 \\ &\Lambda_\gamma(z^2h)=0 \end{aligned}$~~&
  ~~$\begin{aligned} &M(2)\geq0 \\ &\gamma_{01}=\gamma_{10}\\ &\gamma_{02}=\gamma_{11} \\ &a\gamma_{01}+b\gamma_{10}+\gamma_{12}=0 \end{aligned}$~~\\
  \hline
\end{tabular}
\end{center}
\caption{$b$ negative}
\label{tab:table2}
\end{table}

\newpage

\subsection{The TCMP with cubic relation in $M(3)$ of the form $Z^3=itZ+u\overline{Z}$}

 We end this section by considering  another class of cubic column relations in truncated moment problems.

 Set $w=e^{-i\frac{\pi}{4}}z$,  the form $w^3=itw+u\overline{w}$ became $z^3+tz+u\overline{z}=0$.

 We take $t=a$ and $u=b$, then $w^3=itw+u\overline{w}$ owns exactly 7 roots if and only if $u<|t|<2u$, see Table \eqref{tab:table1}.
 Hence if $u\leq 0$ then $\text{rank }M(3)<7$ (more precisely  $\text{rank }M(3)\leq5$); as noted in the introduction this case is not interesting.

In view of Theorems \ref{t.4} and \ref{t.5} we deduce the solution of the TCMP for cubic column relations in $M(3)$ of the form $Z^3=itZ+u\overline{Z}$, where $u<|t|<2u$.

If $0<u<t<2u$; then
\begin{equation*}\begin{aligned}
&h(z,\overline{z})=0,\\
&(z+\overline{z})(z\overline{z}-u)=0,\\
&(e^{\frac{\pi}{4}i}w+e^{-\frac{\pi}{4}i}\overline{w})(\overline{w}w-u)=0,\\
&i(w-i\overline{w})(\overline{w}w-u)=0.
\end{aligned}
\end{equation*}

If $0<u<-t<2u$; then
\begin{equation*}\begin{aligned}
&h(z,\overline{z})=0,\\
&(z-\overline{z})(z\overline{z}-u)=0,\\
&(e^{\frac{\pi}{4}i}w-e^{-\frac{\pi}{4}i}\overline{w})(\overline{w}w-u)=0,\\
&i(w+i\overline{w})(\overline{w}w-u)=0.
\end{aligned}
\end{equation*}

 Now we are able to state the main theorem in \cite{4}.

\begin{theorem}\label{t.6}
Let $\omega\equiv\{\gamma_{ij}\}_{0\leq i+j\leq 6}$, with $\overline{\gamma}_{ij}=\gamma_{ji}$ and $\gamma_{00}>0$, be a truncated complex sequence, and let
 $M(\omega)(3)$ be its associated moment matrix. If $M(\omega)(3)\geq 0$ and has cubic column relation of the form $Z^3=itZ-u\overline{Z}$ , with $t,u\in \mathbb{R}$ and $u<t<2u$.
 Then the following statements are equivalent:
\begin{itemize}
  \item [$i)$] There exists a representing measure for $\omega$.

  \item [$ii)$] There exists a representing measure for the RDIS $\gamma\equiv\{\gamma_{ij}\}_{i,j\geq 0}$, whose initial conditions and characteristic
  polynomial are $\{\gamma_{ij}\}_{0\leq i\leq j\leq 2}$ and $P(z,\overline{z})=z^3-itz-u\overline{z}$, respectively .

  \item [$iii)$]
    \begin{equation*}\begin{cases}
\Lambda_\omega(h)&=0,\\
\Lambda_\omega(zh)&=0.
    \end{cases}\end{equation*}

  \item [$iv)$]
    \begin{equation*}\begin{cases}
\mathfrak{Re}\gamma_{12}-\mathfrak{Im}\gamma_{12}&=u(\mathfrak{Re}\gamma_{01}-\mathfrak{Im}\gamma_{01}),\\
\gamma_{22}&=(t+u)\gamma_{11}-2u\mathfrak{Im}\gamma_{02}.
    \end{cases}\end{equation*}

  \item [$v)$] $h(z,\overline{z})=i(z-i\overline{z})(z\overline{z}-u)\in \mathcal{P}_\gamma$, where $\gamma$ is the RDIS defined in $ii)$.
  \item [$vi)$]  $Z^2 \overline{Z} -iZ \overline{Z}^2 -uZ +iu \overline{Z}=0$.
\end{itemize}
\end{theorem}

\begin{theorem}\label{t.6}
Let $\omega\equiv\{\gamma_{ij}\}_{0\leq i+j\leq 6}$, with $\overline{\gamma}_{ij}=\gamma_{ji}$ and $\gamma_{00}>0$, be a truncated complex sequence, and let
 $M(\omega)(3)$ be its associated moment matrix. If $M(\omega)(3)\geq 0$ and has cubic column relation of the form $Z^3=itZ-u\overline{Z}$ , with $t,u\in \mathbb{R}$ and $u<-t<2u$.
 Then the following statements are equivalent:
\begin{itemize}
  \item [$i)$] There exists a representing measure for $\omega$.

  \item [$ii)$] There exists a representing measure for the RDIS $\gamma\equiv\{\gamma_{ij}\}_{i,j\geq 0}$, whose initial conditions and characteristic
  polynomial are $\{\gamma_{ij}\}_{0\leq i\leq j\leq 2}$ and $P(z,\overline{z})=z^3-itz-u\overline{z}$, respectively .

  \item [$iii)$]
    \begin{equation*}\begin{cases}
\Lambda_\omega(h)&=0,\\
\Lambda_\omega(zh)&=0.
    \end{cases}\end{equation*}

  \item [$iv)$]
    \begin{equation*}\begin{cases}
\mathfrak{Re}\gamma_{12}+\mathfrak{Im}\gamma_{12}&=u(\mathfrak{Re}\gamma_{01}+\mathfrak{Im}\gamma_{01}),\\
\gamma_{22}&=(u-t)\gamma_{11}+2u\mathfrak{Im}\gamma_{02}.
    \end{cases}\end{equation*}

  \item [$v)$] $h(z,\overline{z})=i(z+i\overline{z})(z\overline{z}-u)\in \mathcal{P}_\gamma$, where $\gamma$ is the RDIS defined in $ii)$.
  \item [$vi)$]  $Z^2 \overline{Z} +iZ \overline{Z}^2 -uZ -iu \overline{Z}=0$.
\end{itemize}
\end{theorem}

 \section{ Solving the TCMP with column dependence relations of the form $Z^{k+1}= \sum\limits_{0\leq n+ m \leq k} a_{nm} \overline{Z}^n Z^m$ ($a_{nm} \in \mathbb{C}$)}

 In this section, we involve the RDIS to solve the TCMP associated with the truncated sequence $\gamma \equiv \{ \gamma_{ij}\}_{0\leq i+j\leq 2k+2}$, with
 $\gamma_{ij} = \gamma_{ji}$ and  $M(k+1)(\gamma)$ has a column relation of the form
 \begin{equation}\label{last}
 Z^{k+1}= \sum\limits_{0\leq n+m \leq k} a_{nm} \overline{Z}^n Z^m \phantom{aaa} (a_{nm} \in \mathbb{C}, \text{ for all } n,m\in\mathbb{Z}_+ \text{ and } n+ m \leq k).
  \end{equation}

  According to \eqref{last}, we have $\gamma_{i+k+1,j} = \sum\limits_{0\leq n+m \leq k} a_{nm} \gamma_{n+i, m+j}$, for all $i, j \in \mathbb{Z}_+$ such that $i+j \leq k+1$. Hence $\gamma$ is a subsequence of the RDIS $\widetilde{\gamma}$, defined by the initial conditions $\{ \gamma_{ij}\}_{0\leq i\leq j\leq k}$ and by the characteristic polynomial $p(z, \overline{z}) = z^{k+1} -\sum\limits_{0\leq n+m \leq k} a_{nm} \overline{z}^n z^m$.

  We give now the main result of this section.
  \begin{theorem}\label{th.last}
  Let $M(k+1)(\gamma)$  has a column dependence relations of the form $Z^{k+1}= \sum\limits_{0\leq n+m \leq k} a_{nm} \overline{Z}^n Z^m$ and let $\widetilde{\gamma}$ be a RDIS defined as above. Then $M(k+1)(\gamma)$ admits a representing measure if and only if $M(2k)(\widetilde{\gamma})$ is positive semidefinite.
  \end{theorem}

  {\it proof.}
  Suppose that $\gamma \equiv \{\gamma_{ij}\}_{0\leq i+j\leq 2k+2}$ is a moment sequence, then there exists a positive Borel measure $\mu$ verifies the relation  $$ \gamma_{ij} = \int \overline{z}^i z^j d\mu \phantom{text} (i+j \leq 2k+2).$$
  Set $p(\overline{z}, z)= z^{k+1} -\sum\limits_{0\leq n+m \leq k} a_{nm} \overline{z}^n z^m$. Since $M(k+1)(\gamma)$  has a column dependence relations of the form $p(\overline{Z}, Z) = 0$, then
   $$\int \overline{z}^i z^j p(\overline{z}, z) d\mu = \gamma_{i+k+1, j} -\sum\limits_{0\leq n+m\leq k} a_{nm} \gamma_{i+n, j+m} = 0,$$
  for every $i,j \in \mathbb{Z}_+$ such that $i+j \leq k+1$. Hence
  \begin{equation*}\begin{aligned}
  \int \mid p(\overline{z}, z)\mid^2 d\mu &= \int \overline{p(\overline{z}, z)} p(\overline{z}, z) d\mu \\
  &= \int \overline{z}^{k+1}p(\overline{z}, z) d\mu -\sum\limits_{0\leq n+m\leq k} \overline{a}_{nm} \int \overline{z}^n z^m p(\overline{z}, z) d\mu \\
  &=0,
  \end{aligned}
  \end{equation*}
  thus $supp \mu \subset Z(p)$. It follows, from Proposition \eqref{p.1}, that $\widetilde{\gamma}$ is a moment sequence, and obviously $M(2k)(\widetilde{\gamma})$ is positive semidefinite. Conversely, if $M(2k)(\widetilde{\gamma}) \geq 0$, then Corollary \eqref{c.t.2} yields that $\widetilde{\gamma}$ has a representing measure, and thus $\gamma \equiv \{\gamma_{ij}\}_{0\leq i+j\leq 2k+2}$ is a moment sequence.
  \fin

  On account of Theorem \eqref{th.last} we can formulate the following corollary, which proved a complete solution to  the truncated moment problems with cubic column relations.
  \begin{cor}
  Let $\gamma \equiv \{\gamma_{ij}\}_{0\leq i+j\leq6}$ (with $\overline{\gamma_{ij}} = \gamma_{ji}$) be a complex numbers, let $M(3)(\gamma)$ admits a cubic column relations of the form $Z^3 = \sum\limits_{0 \leq i+j\leq 2} a_{ij}\overline{Z}^i Z^j$ and let $\widetilde{\gamma}$ be the RDIS defined by the initial condition $\{\gamma_{ij}\}_{0\leq i,j\leq2}$ and by the characteristic polynomial $z^3 -\sum\limits_{0 \leq i+j\leq 2} a_{ij}\overline{z}^i z^j$. Then $\gamma$ admits a representing measure if and only if $M(4)(\widetilde{\gamma}) \geq 0$.
  \end{cor}


\begin{thebibliography}{99}
\bibitem{bt}  C. Bayer and J. Teichmann, The proof of Tchakaloff's theorem Proc. Amer. Math. Soc.  134  (2006)  3035-3040.
\bibitem{crz}   C.E. Chidume, M. Rachidi, E.H. Zerouali, Solving the general truncated moment problem by r-generalized Fibonacci sequences methods,                 J. Math. Anal. Appl. 256  (2001) 625-635.
\bibitem{1} R. E. Curto and L. A. Fialkow, solution of the truncated complex moment problem for flat data, Mem. Amer. Math. Soc. 119 (1996).
\bibitem {2} R. E. Curto and L. A. Fialkow, Flat extensions of positive moment matrices: Recursively generated relation, Mem. Amer. Math. Soc. 648 (1998).
\bibitem{3} R. E. Curto and L. Fialkow, Flat extensions of positive moment matrices: Relations in analytic or conjugate terms, Operator Th.: Adv.  Appl. 104 (1998), 59-82.
\bibitem{cy1} R. E. Curto and S. Yoo, Concrete Solution to the Nonsingular Quartic Binary Moment Problem, Proc. Amer. Math. Soc. v.144 (2015)    249-258.
\bibitem{4} R. E. Curto and S. Yoo, Cubic column relations in truncated moment problems, J. Funct. Anal. 266 (2014)  1611-1626.
\bibitem{9} R. E. Curto and L. Fialkow, A Duality Proof of Tchakoloff's Theorem, J.Math. Anal. Appl. 269 (2002) 519-532.
\bibitem{cr1} R. Curto and L. Fialkow, Solution of the truncated hyperbolic moment problem, Integral Equations Operator Theory 52 (2005) 181-218.
\bibitem{cfm} R. Curto, L. Fialkow and H.M. Muller, The extremal truncated moment problem, Integral Equations Operator Theory 60(2008) 177-200.
\bibitem{5} L. Fialkow and J. Nie, Positivity of Riesz functionals and solutions of quadratic and quartic moment problems, J. Funct. Anal. 258(2010) 328-356.
\bibitem{Smu} J.L. Smul'jan, An operator Hellinger integral (Russian), Mat. Sb. 91 (1959) 381-430.
\end{thebibliography}
\end{document}